\documentclass[12pt,reqno,draft]{amsart}
\usepackage{amsmath,amsthm,amssymb,amscd}

\begin{document}

\newcommand{\TITLE}{The Sign of an Elliptic Divisibility Sequence}
\newcommand{\TITLERUNNING}{The Sign of an Elliptic Divisibility Sequence}
\newcommand{\DATE}{February 2004}

%%%%%%%% Set Up Theorem-Style Formats %%%%%%%%%%%%%%
\theoremstyle{plain}
\newtheorem{theorem}{Theorem}
\newtheorem{conjecture}[theorem]{Conjecture}
\newtheorem{proposition}[theorem]{Proposition}
\newtheorem{lemma}[theorem]{Lemma}
\newtheorem{corollary}[theorem]{Corollary}

\theoremstyle{definition}
% A * surpresses numbering, for example
% \newtheorem*{definition}{Definition}
\newtheorem{definition}{Definition}

\theoremstyle{remark}
\newtheorem{remark}{Remark}
\newtheorem{example}{Example}
\newtheorem*{question}{Question}
\newtheorem*{acknowledgement}{Acknowledgements}

\def\BigStrut{\vphantom{$(^{(^(}_{(}$}} % Add space in tables

%%%%%%%% Set Up Environment for Notation %%%%%%%%%%%%%%
% This is currently set to allow quite wide items to be defined
\newenvironment{notation}[0]{%
  \begin{list}%
    {}%
    {\setlength{\itemindent}{0pt}
     \setlength{\labelwidth}{4\parindent}
     \setlength{\labelsep}{\parindent}
     \setlength{\leftmargin}{5\parindent}
     \setlength{\itemsep}{0pt}
     }%
   }%
  {\end{list}}

%%%%%%%% Set Up Environment for Parts in Theorems %%%%%%%%%%%%%%
\newenvironment{parts}[0]{%
  \begin{list}{}%
    {\setlength{\itemindent}{0pt}
     \setlength{\labelwidth}{1.5\parindent}
     \setlength{\labelsep}{.5\parindent}
     \setlength{\leftmargin}{2\parindent}
     \setlength{\itemsep}{0pt}
     }%
   }%
  {\end{list}}
% Use \Part{(a)}, instead of \item[(a)], to ensure upright font
\newcommand{\Part}[1]{\item[\upshape#1]}

%%%%%%%%%%%%%%%%%%
% Greek Alphabet %
%%%%%%%%%%%%%%%%%%
\renewcommand{\a}{\alpha}
\renewcommand{\b}{\beta}
\newcommand{\g}{\gamma}
\renewcommand{\d}{\delta}
\newcommand{\e}{\epsilon}
\newcommand{\f}{\phi}
\newcommand{\fhat}{{\hat\phi}}
\renewcommand{\l}{\lambda}
\renewcommand{\k}{\kappa}
\newcommand{\lhat}{\hat\lambda}
\newcommand{\m}{\mu}
\renewcommand{\o}{\omega}
\renewcommand{\r}{\rho}
\newcommand{\rbar}{{\bar\rho}}
\newcommand{\s}{\sigma}
\newcommand{\sbar}{{\bar\sigma}}
\renewcommand{\t}{\tau}
\newcommand{\z}{\zeta}

\newcommand{\D}{\Delta}
\newcommand{\F}{\Phi}
\newcommand{\G}{\Gamma}

%%%%%%%%%%%%%%%%%%%%
% Fraktur Alphabet %
%%%%%%%%%%%%%%%%%%%%
\newcommand{\ga}{{\mathfrak{a}}}
\newcommand{\gb}{{\mathfrak{b}}}
\newcommand{\gc}{{\mathfrak{c}}}
\newcommand{\gd}{{\mathfrak{d}}}
\newcommand{\gm}{{\mathfrak{m}}}
\newcommand{\gn}{{\mathfrak{n}}}
\newcommand{\gp}{{\mathfrak{p}}}
\newcommand{\gq}{{\mathfrak{q}}}
\newcommand{\gP}{{\mathfrak{P}}}
\newcommand{\gQ}{{\mathfrak{Q}}}

%%%%%%%%%%%%%%%%%%%%%%%%%
% Calligraphic Alphabet %
%%%%%%%%%%%%%%%%%%%%%%%%%
\def\Acal{{\mathcal A}}
\def\Bcal{{\mathcal B}}
\def\Ccal{{\mathcal C}}
\def\Dcal{{\mathcal D}}
\def\Ecal{{\mathcal E}}
\def\Fcal{{\mathcal F}}
\def\Gcal{{\mathcal G}}
\def\Hcal{{\mathcal H}}
\def\Ical{{\mathcal I}}
\def\Jcal{{\mathcal J}}
\def\Kcal{{\mathcal K}}
\def\Lcal{{\mathcal L}}
\def\Mcal{{\mathcal M}}
\def\Ncal{{\mathcal N}}
\def\Ocal{{\mathcal O}}
\def\Pcal{{\mathcal P}}
\def\Qcal{{\mathcal Q}}
\def\Rcal{{\mathcal R}}
\def\Scal{{\mathcal S}}
\def\Tcal{{\mathcal T}}
\def\Ucal{{\mathcal U}}
\def\Vcal{{\mathcal V}}
\def\Wcal{{\mathcal W}}
\def\Xcal{{\mathcal X}}
\def\Ycal{{\mathcal Y}}
\def\Zcal{{\mathcal Z}}

%%%%%%%%%%%%%%%%%%%%%%%%%%%%
% Blackboard Bold Alphabet %
%%%%%%%%%%%%%%%%%%%%%%%%%%%%
\renewcommand{\AA}{\mathbb{A}}
\newcommand{\BB}{\mathbb{B}}
\newcommand{\CC}{\mathbb{C}}
\newcommand{\FF}{\mathbb{F}}
\newcommand{\GG}{\mathbb{G}}
\newcommand{\NN}{\mathbb{N}}
\newcommand{\PP}{\mathbb{P}}
\newcommand{\QQ}{\mathbb{Q}}
\newcommand{\RR}{\mathbb{R}}
\newcommand{\ZZ}{\mathbb{Z}}

%%%%%%%%%%%%%%%%%%%%%%%%%%
% Boldface Math Alphabet %
%%%%%%%%%%%%%%%%%%%%%%%%%%
\def \bfa{{\mathbf a}}
\def \bfb{{\mathbf b}}
\def \bfc{{\mathbf c}}
\def \bfe{{\mathbf e}}
\def \bff{{\mathbf f}}
\def \bfF{{\mathbf F}}
\def \bfg{{\mathbf g}}
\def \bfn{{\mathbf n}}
\def \bfp{{\mathbf p}}
\def \bfr{{\mathbf r}}
\def \bfs{{\mathbf s}}
\def \bft{{\mathbf t}}
\def \bfu{{\mathbf u}}
\def \bfv{{\mathbf v}}
\def \bfw{{\mathbf w}}
\def \bfx{{\mathbf x}}
\def \bfy{{\mathbf y}}
\def \bfz{{\mathbf z}}
\def \bfX{{\mathbf X}}
\def \bfU{{\mathbf U}}
\def \bfmu{{\boldsymbol\mu}}

%%%%%%%%%%%%%%%%%%%%%%%%
% Barred Math Alphabet %
%%%%%%%%%%%%%%%%%%%%%%%%
\newcommand{\Gbar}{{\bar G}}
\newcommand{\Kbar}{{\bar K}}
\newcommand{\Obar}{{\bar O}}
\newcommand{\Pbar}{{\bar P}}
\newcommand{\Qbar}{{\bar Q}}
\newcommand{\QQbar}{{\bar{\QQ}}}

%%%%%%%%%%%%%%%%%%%%%%%%%%%%%%
% Miscellaneous New Commands %
%%%%%%%%%%%%%%%%%%%%%%%%%%%%%%

\newcommand{\Aut}{\operatorname{Aut}}
\newcommand{\Disc}{\operatorname{Disc}}
\renewcommand{\div}{\operatorname{div}}
\newcommand{\Div}{\operatorname{Div}}
\newcommand{\Etilde}{{\tilde E}}
\newcommand{\End}{\operatorname{End}}
\newcommand{\Fix}{\operatorname{Fix}}
\newcommand{\Frob}{\operatorname{Frob}}
\newcommand{\Gal}{\operatorname{Gal}}
\newcommand{\GCD}{{\operatorname{GCD}}}
\renewcommand{\gcd}{{\operatorname{gcd}}}
\newcommand{\hhat}{{\hat h}}
\newcommand{\Hom}{\operatorname{Hom}}
\newcommand{\Ideal}{\operatorname{Ideal}}
\newcommand{\Image}{\operatorname{Image}}
\newcommand{\longhookrightarrow}{\lhook\joinrel\relbar\joinrel\rightarrow}
\newcommand{\LS}[2]{\genfrac(){}{}{#1}{#2}}  % Legendre symbol
\newcommand{\MOD}[1]{~(\textup{mod}~#1)}
\newcommand{\Norm}{\operatorname{N}}
\newcommand{\NS}{\operatorname{NS}}
\newcommand{\notdivide}{\nmid}
\newcommand{\ord}{\operatorname{ord}}
\newcommand{\Parity}{\operatorname{Parity}}
\newcommand{\Pic}{\operatorname{Pic}}
\newcommand{\Proj}{\operatorname{Proj}}
\newcommand{\rank}{\operatorname{rank}}
\newcommand{\res}{\operatornamewithlimits{res}}
\newcommand{\Resultant}{\operatorname{Resultant}}
\renewcommand{\setminus}{\smallsetminus}
\newcommand{\sign}{\operatorname{Sign}}
\newcommand{\Spec}{\operatorname{Spec}}
\newcommand{\Support}{\operatorname{Support}}
\newcommand{\tors}{{\textup{tors}}}
\newcommand\W{W^{\vphantom{1}}} 
\newcommand{\<}{\langle}
\renewcommand{\>}{\rangle}

\hyphenation{para-me-tri-za-tion}

%%%%%%%%%%%%%%%  Topmatter %%%%%%%%%%%%%%%%%%

\title[\TITLERUNNING]{\TITLE}
\date{\DATE}
\author{Joseph H. Silverman}
\address{Mathematics Department, Box 1917, Brown University,
Providence, RI 02912 USA}
\email{jhs@math.brown.edu}
\author{Nelson Stephens}
\address{Department of Mathematics, Royal Holloway, University of London,
Egham, Surrey, TW20 0EX, UK}
\email{nelson.stephens@rhul.ac.uk}

%% \thanks{The first author's research supported by NSA grant ***}

%%%%%%%%%%%%%%%%%%%%%%%%%%%%%%%%%%%%%%%%%%%%%%%%%%%%%%%%%%%%%%%%%%%%%%
%%% Text (non-TeX) Abstract
%% An elliptic divisibility sequence (EDS) is a sequence of integers
%% W_0,W_1,W_2,... generated by the nonlinear recursion satisfied by the
%% division polyomials of an elliptic curve. We give a formula for the
%% sign of W_n for unbounded nonsingular elliptic divisibility sequences.
%% A typical case is Sign(W_n) = (-1)^[n*b] for an irrational real number
%% b, where [x] denotes the greatest integer in x.  As an application, we
%% show that the associated sequence of absolute values
%% |W_1|,|W_2|,|W_3|,... cannot be realized as the
%% sequence counting fixed points of any (abstract) dynamical system.
%%%%%%%%%%%%%%%%%%%%%%%%%%%%%%%%%%%%%%%%%%%%%%%%%%%%%%%%%%%%%%%%%%%%%%

\begin{abstract}
An elliptic divisibility sequence (EDS) is a sequence of
integers~$(W_n)_{n\ge0}$ generated by the nonlinear recursion
satisfied by the division polyomials of an elliptic curve. We give a
formula for the sign of~$W_n$ for unbounded nonsingular EDS, a typical
case being $\sign(W_n)=(-1)^{\lfloor n\b\rfloor}$ for an irrational
number $\b\in\RR$. As an application, we show that the associated
sequence of absolute values~$(|W_n|)$ cannot be realized as the fixed
point counting sequence of any abstract dynamical system.
\end{abstract}

\maketitle

%%%%%%%%%%%%%%%%%%%%%%%%%%%%%%%%%%%%%%%%%%%%%%%%%%%%%%%%%%%%%%%%%%%%%%
\section*{Introduction}
A \textit{divisibility sequence} is a sequence~$(D_n)_{n\ge0}$ of
positive integers with the property that
\[
  m|n \Longrightarrow D_m|D_n.
\]
Classical examples include sequences of the form $a^n-1$ and various
other linear recurrence sequences such as the Fibonacci sequence.
See~\cite{BPvdP} for a complete classification of
divisibility sequences arising from linear recurrences.
\par
There are also natural divisibility sequences associated to nonlinear
recurrence relations. The most famous such relation comes from the
recursion formula for division polynomials on an elliptic curve.

\begin{definition}
An \textit{elliptic divisibility sequence} (abbreviated EDS)
is a divisibility sequence  $W_0,W_1,W_2,\ldots$ that satisfies the formula
\begin{multline}
  \label{equation:EDSrecursion}
  \W_{m+n}\cdot \W_{m-n} = \W_{m+1}\cdot \W_{m-1}\cdot W_n^2
     - \W_{n+1}\cdot \W_{n-1}\cdot W_m^2 \\
  \text{for all $m\ge n\ge 1$.}
\end{multline}
Note that the definition forces $W_0=0$ (put $m=n$) and, except in
some degenerate cases,  $W_1=\pm1$ (put $n=1$). It is also 
possible to use the recursion~\eqref{equation:EDSrecursion} to 
extend the EDS backwards to negative indices, leading to
the identity $W_{-n}=-W_n$.
\end{definition}

The arithmetic properties of elliptic divisibility sequences were
first studied in detail by Morgan Ward~\cite{Ward1,Ward2}.  (See also
\cite{Durst,EGW,EPSW,EW1,EW2,SHIP,SW}.)  Ward describes a
number of degenerate cases for EDS and shows that all other EDS are
associated in a precise way to a pair~$(E,P)$ consisting of an
elliptic curve~$E/\QQ$ and a point~$P\in E(\QQ)$. These nondegenerate
EDS are called \emph{nonsingular}. (See
Definition~\ref{definition:nonsingularityofEDS} for details.)

\begin{example}
\label{example:intro:EDS37}
The simplest (unbounded) nonsingular elliptic divisibility sequence is
the sequence
\begin{multline}
  1, 1, -1, 1, 2, -1, -3, -5, 7, -4, -23, 29, 59, 129, -314, \\
  -65, 1529, -3689, -8209, -16264,\dots
\end{multline}
See Section~\ref{section:numericalexamples} for further examples.
\end{example}

Let~$E/\QQ$ be an elliptic curve given by a Weierstrass equation and
let $P\in E(\QQ)$ be a nontorsion point. For each $n\ge1$, we can
write the $x$-coordinate of~$nP$ in lowest terms as
$x_{nP}=A_{nP}/D_{nP}^2$. It is not hard to prove that $(D_{nP})$ is a
divisibility sequence, and as shown by Shipsey~\cite{SHIP}, it
frequently turns out that~$(D_{nP})$ is
an elliptic divisibility sequence, i.e., it satisfies the
recursion~\eqref{equation:EDSrecursion}, provided that each~$D_{nP}$ 
is chosen with the correct sign. See Section~\ref{section:EDSandEC}
for further details of this connection.
\par
Thus \textit{a priori}, the geometric construction of ``elliptic
divisibility sequences'' via rational points on elliptic curves yields
sequences $(D_{nP})$ of positive integers, while the algebraic
construction via the recursion~\eqref{equation:EDSrecursion} yields
sequences~$(W_n)$ of signed integers.  It is thus of interest to gain
some understanding of how the signs of the terms of an EDS vary.
Our first main result answers this question.

\begin{theorem}
\label{theorem:intro:signofEDS}
Let $(W_n)$ be an unbounded nonsingular elliptic divisibility
sequence. Then possibly after replacing~$(W_n)$ by the related
sequence~$\bigl((-1)^{n-1}W_n\bigr)$, there is an irrational number $\b\in\RR$
so that the sign of~$W_n$ is given by one of the following formulas:
\begin{align*}
  \sign(W_n) &=  (-1)^{\lfloor n\b\rfloor} \quad\text{for all $n$.}\\
  \sign(W_n) &= \begin{cases}
          (-1)^{\lfloor n\b\rfloor+n/2} &\text{if $n$ is even,} \\
          (-1)^{(n-1)/2} &\text{if $n$ is odd.} \\
   \end{cases}
\end{align*}
{\upshape(}Here~$\lfloor t\rfloor$ denotes the greatest integer
in~$t$.{\upshape)}
\end{theorem}

We will prove a more precise theorem in which we use the
parametrization of the associated elliptic curve as a real Lie group
to describe the number~$\b$ and to determine which formula to use.
See Theorem~\ref{theorem:signofEDS} for details.

Our second main result concerns the realizability of an EDS. In
general, a sequence~$(U_n)$ of positive integers is called
\emph{realizable} if there exists a set~$X$ and a function $T:X\to X$
so that
\[
  U_n = \#\bigl\{x\in X : T^n(x)=x\bigr\}
  \qquad\text{for all $n\ge1$.}
\]
Thus realizable sequences are those arising from the theory of (abstract)
dynamical systems. It is clear that an EDS cannot itself be
realizable, since an EDS always has both positive and negative terms.
Our second main result shows that the absolute values of the terms in
an EDS are also not realizable.

\begin{theorem}
\label{theorem:intro:EDSnotrealizable}
Let~$(W_n)$ be an unbounded nonsingular elliptic divisibility
sequence.  Then the associated positive sequence~$(|W_n|)$ is not
realizable.
\end{theorem}

\section{Preliminaries on elliptic divisibility sequences}
\label{section:preliminaries}

A divisibility sequence~$(D_n)$ is called \emph{normalized} if $D_0=0$
and $D_1=1$. Notice that $D_n|D_0$, so if $D_0\ne0$, then~$D_n$ is
bounded.  Further, since~$D_1|D_n$ for all~$n\ge1$, a divisibility
sequence may be normalized by replacing~$D_n$ with~$D_n/D_1$. We will
assume henceforth that all of our elliptic divisibility sequences are
normalized.

An elliptic divisibility sequence~$(W_n)$ is required to satisfy 
the recursion~\eqref{equation:EDSrecursion} for all $m\ge n\ge1$.
It is not hard to show that it suffices that~$(W_n)$ satisfy
the two relations
\begin{align}
  \label{equation:EDSrecursionodd}
  \W_{2n+1} &=\W_{n+2}W_n^3 - \W_{n-1}W_{n+1}^3, \\
  \label{equation:EDSrecursioneven}
  \W_{2n}\W_2 &= \W_n\left(\W_{n+2}W_{n-1}^2 -\W_{n-2}W_{n+1}^2\right).
\end{align}
In particular, an EDS is determined by the values
of~$W_2,W_3,W_4$. Further, a triple~$W_2,W_3,W_4$ with $W_2W_3\ne0$
gives an EDS if and only if~$W_2|W_4$. (See~\cite{Ward1}.)  
\par
We observe that if~$(W_n)$ is an EDS, then $\bigl((-1)^{n-1}W_n\bigr)$
is also an EDS. We call $\bigl((-1)^{n-1}W_n\bigr)$ the \emph{inverse
EDS} to~$(W_n)$, since in Ward's characterization (see
Theorem~\ref{theorem:EDS=sigma}) of EDS via elliptic functions, it
corresponds to replacing the generating point~$z$ by its additive
inverse~$-z$.
\par
The elliptic divisibility sequences that are associated to
elliptic curves are characterized by the following property.

\begin{definition}
\label{definition:nonsingularityofEDS}
The \textit{discriminant} of an elliptic divisibility sequence
sequence~$(W_n)$ is defined by the formula
\begin{align}
  \Disc(W) = {\W_4}W_2^{15} &- {W_3^{3}}W_2^{12} +
  {{3}W_4^{2}}W_2^{10}  - {{{20}\W_4}W_3^{3}}W_2^{7} \notag \\
  &+ {{3}W_4^{3}}W_2^{5} 
    + {{16}W_3^{6}}W_2^{4} 
  + {{{8}W_4^{2}}W_3^{3}}W_2^{2} +  {W_4^{4}}.
  \label{equation:Disc(W)}
\end{align}
The elliptic divisibility sequence~$(W_n)$ is said to be \emph{nonsingular}
if
\[
  W_2\ne0,\qquad W_3\ne0,\qquad\text{and}\qquad \Disc(W)\ne0.
\]
(Cf. Ward~\cite[equation (19.3)]{Ward1}.)
\end{definition}

See \cite{Durst,EGW,EPSW,EW1,EW2,SHIP,SW,Ward1,Ward2} for additional
material on elliptic divisibility sequences.  Ward proves that
nonsingular elliptic divisibility sequences arise as values of the
division polynomials of an elliptic curve. The complicated
expression~\eqref{equation:Disc(W)} defining $\Disc(W)$ is
(essentially) the discriminant of the elliptic curve associated to the
sequence~$(W_n)$. See~\cite{Ward1} for details and formulas, some of
which are recalled in Appendix~\ref{appendix:EDSformulas}. We also
note that Ward gives a complete characterization of all singular EDS.

%%%%%%%%%%%%%%%%%%%%%%%%%%%%%%%%%%%%%%%%%%%%%%%%%%%%%%%%%%%%%%%%%%%%%%
\section{Elliptic divisibility sequences and elliptic functions}

We recall Ward's fundamental result relating EDS to values of elliptic
functions.

\begin{theorem}
\label{theorem:EDS=sigma}
Let $(W_n)$ be a nonsingular elliptic divisibility sequence. Then
there is a lattice $L\subset\CC$ and complex number $z\in\CC$
such that
\[
  W_n = \frac{\sigma(nz,L)}{\sigma(z,L)^{n^2}}
  \qquad\text{for all $n\ge1$,}
\]
where~$\sigma(z,L)$ is the Weierstrass $\s$-function associated to the
lattice~$L$. 
\par
Further, the modular invariants~$g_2(L)$ and~$g_3(L)$ associated to
the lattice~$L$ and the Weierstrass values~$\wp(z,L)$ and~$\wp'(z,L)$
associated to the point~$z$ on the elliptic curve~$\CC/L$ are in the
field~$\QQ(W_2,W_3,W_4)$.  In other
words~$g_2(L),g_3(L),\wp(z,L),\wp'(z,L)$ are all defined over the same
field as the terms of the sequence~$(W_n)$.
\end{theorem}
\begin{proof}
See Ward~\cite[Theorems 12.1 and~19.1]{Ward1}. The rational
expressions for~$g_2$ and~$g_3$ in~$\QQ(W_2,W_3,W_4)$ 
are given by \cite[equations~13.6 and~13.7]{Ward1}, while
the rational expressions for~$\wp(z,L)$ and~$\wp'(z,L)$
are given by \cite[equations~13.5 and~13.1]{Ward1}.
For the convenience of the reader, we have reproduced these 
formulas in Appendix~\ref{appendix:EDSformulas}
\par
Ward also shows~\cite[Theorem~22.1]{Ward1} that up to equivalence, a
singular sequence is either the trivial sequence $W_n=n$ or a Lucas
sequence $W_n=\frac{a^n-b^n}{a-b}$. The latter may be viewed as
arising from the degeneration of the $\sigma$-function to a
trigonometric function.
\end{proof}

%%%%%%%%%%%%%%%%%%%%%%%%%%%%%%%%%%%%%%%%%%%%%%%%%%%%%%%%%%%%%%%%%%%%%%
\section{Elliptic divisibility sequences and elliptic curves}
\label{section:EDSandEC}
There is a natural way to attach a divisibility sequence to any
(nontorsion) rational point on an elliptic curve.
(See~\cite{AEC,ATAEC} for basic properties of elliptic curves.)
Let~$E/\QQ$ be an elliptic curve given by a Weierstrass
equation
\begin{equation}
  \label{equation:WE}
  y^2+a_1xy+a_3y=x^3+a_2x^2+a_4x+a_6.
\end{equation}
A nonzero rational point $P\in E(\QQ)$ can be written in
the form
\[
  P = (x_{P},y_{P}) 
    = \left(\frac{A_{P}}{D_{P}^2},\frac{B_{P}}{D_{P}^3}\right)
  \quad\text{with $\gcd(A_P,D_P)=\gcd(B_P,D_P)=1$.}
\]
Assume now that~$P\in E(\QQ)$ is a nontorsion point. 
It is not hard to show that the sequence
\[
  \left\{ D_{nP} : n=1,2,3,\ldots\right\}
\]
is a divisibility sequence. 
\par
Further, it is often the case that
$(D_{nP})$ is an elliptic divisibility sequence, with one important
caveat, namely we must choose the signs correctly.  
For example, suppose that after moving~$P$ to~$(0,0)$, there is
a Weierstrass equation for the curve~$E$ of the form
\begin{multline}
  \label{equation:WEnonsingpoint}
   E:y^{2}+a_{1}xy+a_{3}y=x^{3}+a_{2}x^{2}+a_{4}x \\
  \text{with $a_{1},a_{2},a_{3},a_{4}\in\ZZ$  and $\gcd (a_{3},a_{4})=1.$}
\end{multline}
There is always some multiple~$mP$ for~$P$ for which this can be
accomplished. More precisely, this is possible if and only if~$P$ has
everywhere nonsingular reduction on the N\'eron model of~$E$
over~$\Spec(\ZZ)$.
\par
Now define a sequence~$W_n$ by setting
\[
  W_1=1,\qquad W_2=a_3, \qquad |W_n|=|D_{nP}|\quad\text{for $n\ge2$},
\]
and choosing the subsequent signs of~$W_n$ by the rule
\[
  \sign(W_{n-2}W_{n}) = -\sign(A_{(n-1)P})
  \qquad\text{for $n\ge3$.}
\]
Shipsey~\cite{SHIP} shows that the sequence $(W_{n})$ is an EDS. 
\par
However, if $P\in E(\QQ)$ is singular modulo some prime, 
then it may not be possible to assign signs to~$\{\pm D_{nP}\}$ in order
to make it into an~ EDS. For example, 
Shipsey~\cite{SHIP} (see also~\cite[\S10.3]{EPSW})
shows that the point~$P=(0,0)$ on the curve $E:y^2+27y=x^3+28x^2+27x$
is not associated to an EDS. The point~$P$ is a
singular point on~$E$ modulo~$3$. The point $P'=(-1,0)=3P\in E(\QQ)$ 
is nonsingular modulo every prime, so~$P'$ does give an~EDS.
\par
Conversely, Shipsey~\cite{SHIP} shows that if $(V_{n})$ is an~EDS,
then there is an ellliptic curve $E/\QQ$ of the
form~\eqref{equation:WEnonsingpoint} and an integer~$k\ge1$ so that
the point~$P=(0,0)\in E(\QQ)$ gives a sequence~$(W_n)$ as described
above (i.e., with $|W_n|=D_{nP}$) such that $W_n=V_{nk}/V_k$ for all
$n\ge1$.  For example, the~EDS beginning $[1, 1, 3, 1,\dots]$ comes
from the point $P=(0,0)$ on the curve $E:y^{2}+27y=x^{3}+28x^{2}+27x$
and the value $k=3.$

\begin{remark}
Divisibility sequences of the form $D_n = a^n-1$ may be viewed as
the values of the $n^{\text{th}}$ division polynomial $X^n-1$ of the
multiplicative group~$\GG_m$. The same is true of divisibility
sequences arising from many other linear recurrences, although one
must either use a twisted version of~$\GG_m$ or work in a quadratic
field. Replacing~$\GG_m$ with an elliptic curve naturally leads
to the theory of elliptic divisibility sequences. More generally, one
can create divisibility sequences associated to rational points on
any algebraic group, see~\cite{SilvermanGCDinFGgps} for details.
\end{remark}

%%%%%%%%%%%%%%%%%%%%%%%%%%%%%%%%%%%%%%%%%%%%%%%%%%%%%%%%%%%%%%%%%%%%%%
\section{The sign of an elliptic divisibility sequence}
To ease notation, for any real number~$x$, we define the parity of~$x$
by
\[
  \sign(x) = (-1)^{\Parity(x)}\qquad\text{with $\Parity(x)\in\ZZ/2\ZZ$.}
\]

\begin{theorem}
\label{theorem:signofEDS}
Let $(W_n)$ be an unbounded nonsingular elliptic divisibility
sequence.  Then possibly after replacing~$(W_n)$ by its inverse
sequence~$\bigl((-1)^{n-1}W_n\bigr)$, there is an irrational
number~$\b\in\RR$ so that the parity of~$W_n$ is given by one of the
following formulas:
\begin{align}
  \Parity(W_n)
  &\equiv\left\lfloor n\b\right\rfloor \pmod{2} \quad\text{for all $n$.}
    \label{equation:paritycaseI} \\
  \Parity(W_n) &\equiv \begin{cases}
      \left\lfloor n\b\right\rfloor
         + \frac{n}{2} \pmod{2} &\text{if $n$ is even,} \\
      \frac{n-1}{2} \pmod{2} &\text{if $n$ is odd.} \\
    \end{cases}
    \label{equation:paritycaseII} 
\end{align}
\par
More precisely, use Theorem~\ref{theorem:EDS=sigma} to associate
to~$(W_n)$ an elliptic curve $E/\RR$ and and point~$P\in E(\RR)$.  Fix
an $\RR$-isomorphism $E(\RR)\cong\RR^*/q^\ZZ$ with $q\in\RR$ and
$|q|<1$. Let~$P\leftrightarrow u$, where~$u$ is normalized to satisfy
$q^2<u<1$ if $q<0$ and $|q|<|u|<1$ otherwise.  Then the formula for
the parity of~$W_n$ and the value of~$\b$ are given in the following
table:
\begin{center}
\begin{tabular}{|c|c|c|c|c|c|} \hline
$q$ & $u$ & Formula & $\b$ & $E(\RR)$ & P \\ \hline\hline
$q>0$ & $u>0$ & \eqref{equation:paritycaseI} & $\log_q(u)$ 
    & Disconnected & Identity  \\ 
  & & & & & component \\\hline
$q>0$ & $u<0$ & \eqref{equation:paritycaseII} & $\log_q(|u|)$ 
    & Disconnected & Nonidentity \\ 
  & & & & & component \\\hline
$q<0$ &  & \eqref{equation:paritycaseI} & $\frac{1}{2}\log_{|q|}(u)$ 
    & Connected &  \\ \hline
\end{tabular}
\end{center}
\end{theorem}

\begin{remark}
Theorem~\ref{theorem:signofEDS} is true for any unbounded sequence of
\emph{real numbers}~$(W_n)$ satisfying the
recursion~\eqref{equation:EDSrecursion}, subject to the 
nonsingularity condition~\eqref{definition:nonsingularityofEDS}.
\end{remark}

\begin{proof}[Proof of Theorem~\ref{theorem:signofEDS}]
The proof of Theorem~\ref{theorem:signofEDS} relies on Ward's
fundamental result (Theorem~\ref{theorem:EDS=sigma}) relating ellipitic
divisibility sequences to the values of elliptic functions.
Let $(W_n)$ be an unbounded nonsingular elliptic divisibility sequence,
and choose a lattice~$L$ and complex number~$z$ as in
Theorem~\ref{theorem:EDS=sigma} so that
\begin{equation}
  \label{equation:Wn=sigma}
  W_n = \frac{\sigma(nz,L)}{\sigma(z,L)^{n^2}}.
\end{equation}
Let~$E$ be the associated elliptic curve
\[
  E : Y^2 = 4X^3 - g_2(L)X - g_3(L)
\]
and $P=\bigl(\wp(z,L),\wp'(z,L)\bigr)$ the associated point on~$E$.
Theorem~\ref{theorem:EDS=sigma} tells us that~$E$ is defined
over~$\QQ$ (i.e., $g_2(L),g_3(L)\in\QQ$) and that $P\in E(\QQ)$. In
particular,~$E$ and~$P$ are defined
over~$\RR$. Then~\cite[V.2.3(b)]{ATAEC} says that there exists a
(unique) $q\in\RR$ with $0<|q|<1$ such that there is an
$\RR$-isomorphism
\[
  \psi : \frac{\CC^*}{q^\ZZ} \;\buildrel\sim\over\longrightarrow\; E(\CC).
\]
(See \cite[V.1.1]{ATAEC} for explicit power series defining~$\psi$.)
The fact that~$\psi$ is defined over~$\RR$ means that it gives an
isomorphism \text{$\psi:\RR^*/q^\ZZ\to E(\RR)$}, so the fact that $P\in
E(\RR)$ implies that $\psi^{-1}(P)\in\RR^*/q^\ZZ$.  Let~$u\in\RR^*$ be
a representative for $\psi^{-1}(P)$.
\par
Write $u=e^{2\pi i\a}$ with (say) $\a\in i\RR$ if $u>0$ and
$\a\in\frac{1}{2}+i\RR$ if $u<0$.  Then the~$\s$-function
on~$\CC^*/q^\ZZ$ is given by the formula~\cite[V.1.3]{ATAEC}
\begin{multline}
  \label{equation:sigma=prod}
  \s(u,q)=-\frac{1}{2\pi i}e^{\frac{1}{2}\eta \a^2-\pi i\a}\theta(u,q)
  \quad\text{with}\\
  \theta(u,q)=(1-u)\prod_{m\ge1}
    \frac{(1-q^mu)(1-q^mu^{-1})}{(1-q^m)^2}.
\end{multline}
Here~$\eta$ is a quasiperiod, but its value will not concern us, since
it disappears when we substitute~\eqref{equation:sigma=prod}
into~\eqref{equation:Wn=sigma}. However, it is important to observe
that the $\s$-function in formula~\eqref{equation:Wn=sigma} and the
$\s$-function defined by the formula~\eqref{equation:sigma=prod} may
only be constant multiples of one another, since~$\s(z,L)$ has weight
one (i.e., $\s(cz,cL)=c\s(z,L)$). Hence when we
substitute~\eqref{equation:sigma=prod} into~\eqref{equation:Wn=sigma},
we obtain
\begin{equation}
  \label{equation:Wn=prod}
  W_n = \g^{n^2-1}u^{(n^2-n)/2}
         \frac{\theta(u^n,q)}{\theta(u,q)^{n^2}}
\end{equation}
for some~$\g\in\CC^*$. However, since~$u$,~$q$, and~$W_n$ are all
in~$\RR$, taking~$n=2$ and~$n=3$ shows that~$\g^3\in\RR$
and~$\g^8\in\RR$, so~$\g\in\RR$.
\par
We observe that since $n^2-1\equiv n-1\pmod{2}$ for all $n\in\ZZ$, the
effect of a negative~$\g$ is simply to replace an EDS starting
$[1,a,b,c]$ with the asociated inverse sequence starting $[1,-a,b,-c]$. 
Hence without loss of generality, we may assume that $\g>0$.
\par
Since $|q|<1$, we see that $1-q^m>0$ for all $m\ge1$. Thus in
computing the sign of~$W_n$ using~\eqref{equation:Wn=prod}, we may
discard the $(1-q^m)^2$ factors appearing in the product
expansion~\eqref{equation:sigma=prod} for~$\theta(u,q)$.
We now consider several cases, depending on the sign of~$q$ and~$u$.

\paragraph{\textbf{Case I}: $1>q>0$ and $u>0$}
\leavevmode\newline 
Geometrically, this is the case that~$\RR^*/q^\ZZ$ has two components
and the point~$P=\psi(u)$ is on the identity component. 
The value of the righthand side of~\eqref{equation:Wn=prod} is
invariant under~$u\to q^{\pm1}u$, so we may choose~$u$
to satisfy $q<u<1$. Then 
\[
  1-u>0
  \qquad\text{and}\qquad 1-q^mu^{\pm1} > 0
  \quad\text{for all $m\ge1$,}
\]
so $\theta(u,q)\ge 0$. We next do a similar analysis
for~$\theta(u^n,q)$. The assumption $0<q<u<1$
and the fact that we take $n\ge1$ implies that
\[
  1-q^mu^n > 0
  \qquad\text{for all $m\ge1$.}
\]
So the only sign ambiguity in $\theta(u^n,q)$ comes from the factors
of the form \text{$1-q^mu^{-n}$}. We have
\[
  1-q^mu^{-n}<0 \iff u^n<q^m
  \iff   n\log_q(u) > m.
\]
Hence there are $\lfloor n\log_q(u)\rfloor$ negative signs. This proves
that
\[
  \Parity(W_n) \equiv
  \Parity(\theta(u^n,q)) \equiv \lfloor n\log_q(u)\rfloor
  \pmod{2},
\]
which is the desired result~\eqref{equation:paritycaseI} with
the explicit value $\b=\log_q(u)$.

\paragraph{\textbf{Case II}: $1>q>0$ and $u<0$}
\leavevmode\newline 
Geometrically, we are again in the case that~$\RR^*/q^\ZZ$ has two
components, but now the point~$u$ is on the nonidentity component.  We
may choose~$u$ to satisfy $q<|u|<1$, and then as in Case~I, we see
that~$\theta(u,q)>0$ and that all factors \text{$1-q^mu^n$} are positive.
Further, since $u<0$ and $q>0$, it is clear that
\text{$1-q^mu^{-n}>0$} for all odd values of~$n$. Thus if~$n$ is odd,
we also have $\theta(u^n,q)>0$.
\par
Suppose now that~$n$ is even. Then
\[
  1-q^mu^{-n}<0 \iff |u|^n<q^m
  \iff   n\log_q(|u|) > m.
\]
Hence there are $\lfloor n\log_q(|u|)\rfloor$ negative signs, so
\[
  \Parity(\theta(u^n,q)) \equiv \lfloor n\log_q(|u|)\rfloor \pmod{2}
  \qquad\text{when $n$ is even.}
\]
Finally, since~$u<0$, we observe that
\[
  \Parity\left(u^{(n^2-n)/2}\right)\equiv \frac{n^2-n}{2}
  \equiv \begin{cases}
     n/2\pmod{2} &\text{if $n$ is even,} \\
     (n-1)/2\pmod{2} &\text{if $n$ is odd.} \\
  \end{cases}
\]
Combining these results and substituting into~\eqref{equation:Wn=prod}
yields the desired result~\eqref{equation:paritycaseII} with
the explicit value $\b=\log_q(u)$.

\paragraph{\textbf{Case III}: $q<0$}
\leavevmode\newline
Geometrically, this is the case that the curve~$\RR^*/q^\ZZ$ is connected.
Replacing~$u$ by~$q^ku$ for an appropriate~$k\in\ZZ$, we may assume that
\[
  u > 0 \qquad\text{and}\qquad q^2 < u < 1.
\]
Consider first the factors $1-q^mu^{\pm1}$ of~$\theta(u,q)$. For our
choice of~$u$, we have
\[
  1-q^mu^{\pm1} > 1-|q|^{m-2},
\]
so \text{$1-q^mu^{\pm1}$} is positive except possibly when~$m=1$.
However, when $m=1$, it is also positive, since~$q<0$
and~$u>0$. Hence~$\theta(u,q)>0$.
\par
Next we look at the factors of~$\theta(u^n,q)$.  The factors
\text{$1-q^mu^n$} are positive, since $|q|<1$ and $u<1$.
Further, the factors \text{$1-q^mu^{-n}$} with~$m$ odd are also positive,
since~$q<0$ and~$u>0$. Suppose now that~$m$ is even. Then
\[
  1-q^mu^{-n} < 0
  \iff u^n < |q|^m
  \iff n\log_{|q|}(u) > m.
\]
Thus we get one negative factor in~$\theta(u^n,q)$ for
each even integer smaller than~$n\log_{|q|}(u)$, so
\[
  \Parity(W_n) \equiv \Parity(\theta(u^n,q))
  \equiv \left\lfloor\frac{1}{2}\log_{|q|}(u)\right\rfloor \pmod{2}.
\]
This yields the desired result~\eqref{equation:paritycaseI} with
the explicit value $\b=\frac{1}{2}\log_{|q|}(u)$, which
completes the proof of Theorem~\ref{theorem:signofEDS}.
\end{proof}

%%%%%%%%%%%%%%%%%%%%%%%%%%%%%%%%%%%%%%%%%%%%%%%%%%%%%%%%%%%%%%%%%%%%%%
\section{Numerical examples}
\label{section:numericalexamples}
We give some examples of elliptic divisibility sequences that
illustrate the various cases of Theorem~\ref{theorem:signofEDS}.

\begin{example}
\label{example:EDS37}
The elliptic divisibility sequence starting~$[1,1,-1,1]$ is
\begin{multline}
  \label{equation:EDScond37}
  1, 1, -1, 1, 2, -1, -3, -5, 7, -4, -23, 29, 59, 129, -314, \\
  -65, 1529, -3689, -8209, -16264,\dots
\end{multline}
This is the elliptic divisibility sequence associated to the elliptic
curve $E:y^2+y=x^3-x$ of conductor~37 and the point $P=(0,0)\in
E(\QQ)$. The point~$P$ has canonical height~$\hhat(P)\approx0.0256$.
Among elliptic curves over~$\QQ$ of positive rank, this is the one of
smallest conductor. So in some sense~\eqref{equation:EDScond37} is the
``smallest'' or ``simplest'' possible elliptic divisibility
sequence. The authors of~\cite{EPSW} liken it to the simplest linear
divisibility sequence \text{$2^n-1$}.
\par
There is an isomorphism $E(\RR)\cong\RR^*/q^\ZZ$ with $P\leftrightarrow u$,
where
\[
  q=0.0004654203923\dots\qquad\text{and}\qquad u = -0.09230562888\dots. 
\]
Thus~$E(\RR)$ is disconnected and~$P$ is on
the nonidentity component, so Theorem~\ref{theorem:signofEDS} says
that
\begin{multline*}
  \Parity(W_n) \equiv \begin{cases}
     \lfloor n\b\rfloor +\frac{n}{2} + 1&\text{if $n$ is even,}\\
     \frac{n-1}{2}&\text{if $n$ is odd,}\\
  \end{cases} \\
  \text{with $\b=\log_q(|u|)=0.310541358720\dots$.}
\end{multline*}
Notice that we have added one for the even values of~$n$, since
Theorem~\ref{theorem:signofEDS}  gives the sign of either the
sequence~$(W_n)$ or the sequence $((-1)^{n-1}W_n)$.
\end{example}

\begin{example}
\label{example:EDS43}
The elliptic divisibility sequence starting~$[1,1,1,-1]$ is
\begin{multline}
  \label{equation:EDScond43}
  1, 1, 1, -1, -2, -3, -1, 7, 11, 20, -19, -87, -191, -197, 1018, \\
  2681, 8191, -5841, -81289, -261080,\dots
\end{multline}
This is the elliptic divisibility sequence associated to the elliptic
curve $E:y^2+y=x^3+x^2$ of conductor~43 and the point $P=(0,0)\in
E(\QQ)$. The point~$P$ has canonical height~$\hhat(P)\approx0.0314$.
\par
There is an isomorphism $E(\RR)\cong\RR^*/q^\ZZ$ with $P\leftrightarrow u$,
where
\[
  q=-0.001833413287\dots\qquad\text{and}\qquad u = 0.02931619135\dots. 
\]
Thus~$E(\RR)$ is connected, so Theorem~\ref{theorem:signofEDS} says
that
\[
  \Parity(W_n) \equiv \lfloor n\b\rfloor
  \qquad\text{with $\b=\frac{1}{2}\log_q(|u|)= 0.2800581462\dots$.}
\]
\end{example}

\begin{example}
\label{example:EDS58}
The elliptic divisibility sequence starting~$[1,1,2,1]$ is
\begin{multline}
  \label{equation:EDScond58}
  1, 1, 2, 1, -7, -16, -57, -113, 670, 3983, 23647, 140576, -833503, \\
  -14871471, -147165662, -2273917871, 11396432249, 808162720720, \\
  14252325989831, 503020937289311,\dots
\end{multline}
This is the elliptic divisibility sequence associated to the elliptic
curve $E:y^2+xy=x^3-x^2-x+1$ of conductor~58 and the point $P=(1,0)\in
E(\QQ)$. The point~$P$ has canonical height~$\hhat(P)\approx0.0848$.
The reason that~\eqref{equation:EDScond58} grows so much faster
than~\eqref{equation:EDScond37} or~\eqref{equation:EDScond43} is due
to the larger height of~$P$. One can show that (roughly)
$|W_n|\approx\exp\left(\frac{1}{2}n^2\hhat(P)\right)$.
\par
There is an isomorphism $E(\RR)\cong\RR^*/q^\ZZ$ with $P\leftrightarrow u$,
where
\[
  q=-0.0004429838967\dots\qquad\text{and}\qquad u = 0.02529988312\dots. 
\]
Thus as in the previous example,~$E(\RR)$ is connected and
\[
  \Parity(W_n) \equiv \lfloor n\b\rfloor 
  \qquad\text{with $\b=\frac{1}{2}\log_q(|u|)=0.2380838117\dots$.}
\]
\end{example}

\begin{example}
\label{example:EDS61}
The elliptic divisibility sequence starting~$[1,1,1,2]$ is
\begin{multline}
  \label{equation:EDScond61}
  1, 1, 1, 2, 1, -3, -7, -8, -25, -37, 47, 318, 559, 2023, 7039, -496,\\
  -90431, -314775, -1139599, -8007614,\dots
\end{multline}
This is the elliptic divisibility sequence associated to the elliptic
curve $E:y^2+xy=x^3-2x+1$ of conductor~61 and the point $P=(1,0)\in
E(\QQ)$. The point~$P$ has canonical height~$\hhat(P)\approx0.0396$.
We have $q=-0.00006372107969$ and $u=0.02660268122$, so again~$E(\RR)$
is connected and
\[
  \Parity(W_n) \equiv \lfloor n\b\rfloor 
  \qquad\text{with $\b=\frac{1}{2}\log_q(|u|)=0.1877002949\dots$.}
\]
\end{example}

\begin{example}
\label{example:EDS710}
The elliptic divisibility sequence starting~$[1,2,2,-2]$ is
\begin{multline}
  \label{equation:EDScond710}
  1, 2, 2, -2, -24, -100, -176, 1552, 28448, 248448, 433024, \\
  -47795200, -1682842624, -30121422848, 218738737152,\dots
\end{multline}
This is the elliptic divisibility sequence associated to the elliptic
curve $E:y^2+xy+y=x^3+x^2-416x+3009$ of conductor~710 and the point
$P=(21,53)\in E(\QQ)$. The point~$P$ has canonical
height~$\hhat(P)\approx0.08372$. However, if we write $x(nP)=A_n/D_n^2$,
then it is not true that $|W_n|=D_n$. The difficulty arises because~$P$
is singular modulo~$2$. One can check that $W_n/D_n=\pm2^k$  for all~$n$.
\par
There is an isomorphism $E(\RR)\cong\RR^*/q^\ZZ$ with $P\leftrightarrow u$,
where
\[
  q=0.00002987174044\dots\qquad\text{and}\qquad u = 0.0004951251683\dots. 
\]
Thus~$E(\RR)$ is disconnected and~$P$ is on
the identity component, so Theorem~\ref{theorem:signofEDS} says  
that
\[
  \Parity(W_n) \equiv \lfloor n\b\rfloor+(n-1)
  \qquad\text{with $\b=\log_q(u)=0.7304917812\dots$.}
\]
Notice that as in Example~\ref{example:EDS37}, the theorem gives
the sign of~$(-1)^{n-1}W_n$, so we needed to adjust by \text{$n-1$}.
\end{example}

%%%%%%%%%%%%%%%%%%%%%%%%%%%%%%%%%%%%%%%%%%%%%%%%%%%%%%%%%%%%%%%%%%%%%%
\section{Elliptic divisibility sequences and realizability}
A sequence is said to be realizable if it gives the number of fixed
points of the iterates of some map. As an application of
Theorem~\ref{theorem:signofEDS}, we show in this section that the
absolute value of an elliptic divisibility sequence is not realizable.
See~\cite[Section~11.2]{EPSW} and~\cite{EPW} for further information
about realizability.

\begin{definition}
Let $(U_n)_{n\ge1}$ be a sequence of nonnegative integers. We say
that~$(U_n)$ is \emph{realizable} if there exists a set~$X$ and a
function \text{$\f:X\to X$} with the property that
\[
  U_n = \#\{x\in X : T^n(x)=x\} 
  \qquad\text{for all $n\ge1$.}
\]
\end{definition}

\begin{remark}
It is clear from Theorem~\ref{theorem:signofEDS} that an EDS
cannot be realizable, since~$(W_n)$ will always contain terms that are
negative. It turns out that the sequence of absolute values~$(|W_n|)$
is also not realizable, but the proof, which we give below, is less
direct.
\end{remark}

\begin{remark}
There are many \emph{linear} recurrence sequences that are known to be
realizable. A simple example is the Lucas sequence $L_{n+2}=L_{n+1}+L_n$ with
initial terms $L_1=1$ and $L_2=3$. Linear recurrence sequences grow
linearly exponentially, i.e., $\log|L_n|=O(n)$. Elliptic divisibility
sequences generally grow much more rapidly, namely $\log|W_n|=O(n^2)$. At
the present time, there are no realizable sequences known
that grow quadratic exponentially. Elliptic divisibility sequences, with
their added structure, are thus a natural place to search for rapidly growing
realizable sequences. The main result in this section shows that unfortunately
the EDS recursion never gives a realizable sequence.
\end{remark}

There is an elementary combinatorial characterization of realizable
sequences. We recall that the Dirichlet product of two 
arithmetic functions~$f$ and~$g$ is defined by the formula
\[
  (f*g)(n) = \sum_{d|n} f(d)g\left(\frac{n}{d}\right).
\]

\begin{lemma}
\label{lemma:realizabilitycriterion}
A sequence $(U_n)$ of positive integers is realizable if and only 
\begin{alignat}{2}
  (U*\mu)(n)&\ge0 &\qquad&\text{for all $n\ge1$, and}\\
  (U*\mu)(n)&\equiv0\pmod{n}&&\text{for all $n\ge1$.}
  \label{equation:RScongruence}
\end{alignat}
(Here $\mu(n)$ is the M\"obius function.)
\end{lemma}
\begin{proof}
See~\cite[Lemma~11.3]{EPSW}.
\end{proof}

\begin{remark}
Applying~\eqref{equation:RScongruence} to a prime power~$p^k$ shows
that a realizable sequence necessarily satisfies
\begin{equation}
  \label{equation:RSpadiccong}
  U_{p^k} \cong U_{p^{k-1}} \pmod{p^k}
  \qquad\text{for all $k\ge1$.}
\end{equation}
In particular, the limit
\[
  \lim_{k\to\infty} U_{p^k} \in \ZZ_p\qquad\text{exists in~$\ZZ_p$.}
\]
It would interesting to give a dynamical interpretation
to this limit. For the proof of the nonrealizibility of an EDS,
we will only require the very special case of~\eqref{equation:RSpadiccong}
which asserts that $U_{2^k}\equiv U_{2^{k-1}}\pmod{4}$ for all $k\ge2$.
\end{remark}

The following lemma provides a (weak) EDS counterpart to the $p$-adic
coherence~\eqref{equation:RSpadiccong} of a realizable sequence, at
least for the prime~$p=2$. 

\begin{lemma}
\label{lemma:EDSperiodicmod4}
Let $(W_n)$ be an unbounded EDS, and assume that~$W_2$ and~$W_4$ are
odd.  Then~$W_{2^k}$ is odd for all $k\ge0$, and for any fixed
power~$2^e$, the sequence
$W_{2^k}\bmod{2^e}$ is eventually periodic, i.e., there are integers
$r,K>0$ so that
\[
  W_{2^{r+k}}\equiv W_{2^k} \pmod{2^e}
  \qquad\text{for all $k\ge K$.}
\]
\end{lemma}

\begin{remark}
Lemma~\ref{lemma:EDSperiodicmod4} says that the
sequence~$(W_{2^k}\bmod 2^e)_{k\ge0}$ is eventually periodic, and an
adaptation of the proof of Lemma~\ref{lemma:EDSperiodicmod4} gives a
similar result with~$2$ replaced by other primes.  This suggests that
something stronger might be true.
\begin{question}
Let~$(W_n)$ be an unbounded EDS. Does there exist an exponent~$e\ge1$
so that for each $0\le i<e$, the limit
\begin{equation}
  \label{equation:padiclimit}
  \lim_{k\to\infty} W_{p^{ke+i}}\quad\text{exists in $\ZZ_p$?}
\end{equation}
\end{question}
\noindent
If the underlying elliptic curve~$E$ has split multiplicative
reduction at~$p$, so $E(\QQ_p)\cong\QQ_p^*/q^\ZZ$, then the
limit~\eqref{equation:padiclimit} exists and can be expressed in terms
of Tate's $p$-adic sigma function~(\cite[V~\S3]{ATAEC}) and the
Teichm\"uller character. In general, we conjecture that the
limit~\eqref{equation:padiclimit}  contains interesting $p$-adic
information related to the elliptic Teichm\"uller function
$E(\QQ_p)\to E(\QQ_p)_\tors$ and to the Mazur-Tate $p$-adic sigma
function~\cite{MT}.
\end{remark}

\begin{proof}[Proof (of Lemma~\ref{lemma:EDSperiodicmod4})]
The EDS recurrence~\eqref{equation:EDSrecursion}
allows one to compute an entire EDS from its first four terms. 
More generally, it is possible to generate the subsequence~$(W_{dn})_{n\ge1}$
from the initial four values $W_d,W_{2d},W_{3d},W_{4d}$ using the 
recurrence
\[
  W_d^2 \W_{(n+2)d} \W_{(n-2)d} = W_{2d}^2 \W_{(n+1)d} \W_{(n-1)d}
                - \W_d \W_{3d} W_{nd}^2.
\]
In particular, the following formulas express~$W_{6d}$ and~$W_{8d}$ in
terms of~$W_d,W_{2d},W_{3d},W_{4d}$:
\begin{align}
  W_{6d} &= W_{3d}\left(\frac{W_{2d}^4\W_{4d}}{W_d^5}
           - \frac{\W_{2d}W_{3d}^3}{W_d^4}
           - \frac{W_{4d}^2}{\W_d\W_{2d}} \right) 
    \label{equation:W6drecursion}\\
  W_{8d} &= W_{4d}\left(
  -\frac{2W_{3d}^6}{W_d^6} + \frac{3W_{2d}^3W_{3d}^3\W_{4d}}{W_d^7}
  -\frac{W_{3d}^3W_{4d}^2}{W_d^3W_{2d}^2}
  -\frac{W_{2d}^6W_{4d}^2}{W_d^8} \right)  
    \label{equation:W8drecursion}
\end{align}
\par
We use induction to prove that~$W_{2^k}$ is odd. We are given
that~$W_1$,~$W_2$, and~$W_4$ are odd. Assume now
that~$W_{2^k}$,~$W_{2^{k+1}}$ and~$W_{2^{k+2}}$ are odd for
some~$k\ge0$. We apply~\eqref{equation:W8drecursion} with
$d=2^k$. Since $W_d\equiv W_{2d}\equiv W_{4d}\equiv 1\pmod{2}$ by
hypothesis, we can reduce~\eqref{equation:W8drecursion} modulo~$2$
to obtain
\[
  W_{2^{k+3}} = W_{8d} \equiv 1\cdot(-2+3W_{3d}^3-W_{3d}^3-1)
  \equiv 1 \pmod{2}.
\]
Hence~$W_{2^{k+3}}$ is odd, which completes the induction proof.
\par
Now consider the map
\[
  F:\NN \longrightarrow \left(\frac{\ZZ}{2^e\ZZ}\right)^4,\qquad
  k \longrightarrow 
    \bigl( W_{2^k}, W_{2\cdot2^k}, W_{3\cdot2^k},
       W_{2^e\cdot2^k}\bigr) \bmod 2^e.
\]
We observe that the first two coordinates of~$F(k+1)$ are simply the
second and fourth coordinates of~$F(k)$, while the last two
coordinates of~$F(k+1)$ can be computed from the coordinates of~$F(k)$
using the formulas~\eqref{equation:W6drecursion}
and~\eqref{equation:W8drecursion}. (Note that we are using the fact
that all~$W_{2^k}$ are odd, since otherwise we could not simply
reduce~\eqref{equation:W6drecursion} and~\eqref{equation:W8drecursion}
modulo~$2^e$.) 
\par
We have shown that~$F(k+1)$ is completely determined by~$F(k)$. Since
the range of~$F$ is a finite set, it must eventually repeat a value,
say $F(r+K)=F(K)$ for some $K\ge0$ and $r\ge1$. It follows that
$F(r+k)=F(k)$ for all $k\ge K$.
\end{proof}

\begin{theorem}
\label{theorem:EDSnotrealizable}
Let~$(W_n)$ be an unbounded nonsingular elliptic divisibility sequence.
Then the sequence of absolute values~$(|W_n|)$ is not realizable.
\end{theorem}
\begin{proof}
Suppose that~$(|W_n|)$ is realizable. Applying
Lemma~\ref{lemma:realizabilitycriterion} with $n=2^k$ (i.e.,
equation~\eqref{equation:RSpadiccong} with $p=2$), we see that
\begin{equation}
  \label{equation:2adicstability}
  |W_{2^k}| \equiv |W_{2^{k-1}}|  \pmod{2^k}
  \qquad\text{for all $k\ge1$.}
\end{equation}
\par
In particular, since~$W_1=1$ is odd, we see that~$W_{2^k}$ is odd for
all~$k$. This allows us to apply Lemma~\ref{lemma:EDSperiodicmod4}, 
which we do with~$e=2$. Thus we find integers~$r,K>0$ so that
\begin{equation}
  \label{equation:mod4stability}
  W_{2^{r+k}} \equiv W_{2^k} \pmod{4}
  \qquad\text{for all $k\ge K$.}
\end{equation}
On the other hand,~\eqref{equation:2adicstability} certainly implies that
\begin{equation}
  \label{equation:mod4stabilityabs}
  |W_{2^{r+k}}| \equiv |W_{2^k}| \pmod{4}
  \qquad\text{for all $k\ge1$.}
\end{equation}
The two congruences~\eqref{equation:mod4stability}
and~\eqref{equation:mod4stabilityabs} imply that for
each~$j$ with $0\le j<r$, the quantity
\[
  \text{$\sign(W_{2^{ri+j}})$ is constant for all $i\ge K$.}
\]
\par
However, Theorem~\ref{theorem:signofEDS} tells us that there is
an irrational number~$\b\in\RR$ with the property that
\[
  \sign(W_{2^{k}}) = (-1)^{\lfloor 2^{k}\b\rfloor}
  \qquad\text{for all $k\ge1$,}
\]
so the assumption that~$(|W_n|)$ is realizable leads to the conclusion
that for each $0\le j<r$, the parity of $\lfloor 2^{ri+j}\b\rfloor$ is
constant for all $i\ge K$. In other words, the parity of $\lfloor
2^{k}\b\rfloor$ is eventually periodic.  Then an elemenatry argument
(see Lemma~\ref{lemma:constantparity} below) implies that~$\b\in\QQ$,
contradicting the fact that~$\b$ is irrational. This completes the
proof that~$(|W_n|)$ is not realizable.
\end{proof}

\begin{lemma}
\label{lemma:constantparity}
Let $\b\in\RR$ and suppose that the parity of $\lfloor 2^k\b\rfloor$
is eventually periodic, i.e., there are integers~$r\ge1$ and~$K$ so
that for each $0\le j<r$, the quantity $\lfloor 2^{ri+j}\b\rfloor$ is
constant for all $i\ge K$. Then~$\b\in\QQ$.
\end{lemma}
\begin{proof}
Write the binary expansion of the fractional part of~$\b$ as
\[
   \b = \lfloor\b\rfloor + \sum_{k=1}^\infty \frac{b_k}{2^k}
  \qquad\text{with $b_k \in \{0,1\}$.}
\]
The parity of~$\lfloor 2^k\b\rfloor$ is simply the coefficient~$b_k$,
so our assumption says that the sequence of binary
coefficients~$(b_k)$ is eventually periodic.  Hence~$\b$ is rational.
\end{proof}

\begin{acknowledgement}
The authors would like to thank Graham Everest, Alf van der Poorten,
Igor Shparlinski, Thomas Ward and Gary Walsh for helpful disucssions
about (elliptic) divisibility sequences.
\end{acknowledgement}

%%%%%%%%%%%%%%%%%%%%%%%%%%%%%%%%%%%%%%%%%%%%%%%%%%%%%%%%%%%%%%%%%%%%%%

%%%%%%%%%%%%%%%%%%%%%%%%%%%%%%%%%%%%%%%%%%%%%%%%%%%%%%%%%%%%%%%%%%%%%%
%% \newpage

\appendix

\section{The elliptic curve associated to an elliptic divisiblity sequence}
\label{appendix:EDSformulas}

The material in this appendix is essentially taken from Ward's
paper~\cite{Ward1}.  Let~$(W_n)$ be an elliptic divisibility sequence
starting
\[
  W_1=1,\quad W_2=\a,\quad W_3=\b,\quad W_4=\a\g.
\]
(Notice that we have built in the divisibility of~$W_4$ by~$W_2$.)
Define quantities
\begin{align}
  A &= 3^3\cdot\bigl(
  -\a^{16} - {{4}\g}\a^{12} + \left({16}\b^{3} - {{6}\g^{2}}\right) \a^{8}
  + \left({{8}\g}\b^{3} - {{4}\g^{3}}\right) \a^{4} \notag \\
  &\qquad\qquad\qquad\qquad
          - \left({16}\b^{6} + {{8}\g^{2}}\b^{3} + {\g^{4}}\right) \bigr) 
  \label{equation:formulaforA}\\ 
  %%%%%%%%%%%%%%%%%%%%
  B &= 2\cdot 3^3\cdot\bigl(
  \a^{24} + {{6}\g}\a^{20} - \left({24}\b^{3} - {{15}\g^{2}}\right) \a^{16}
  - \left({{60}\g}\b^{3} - {{20}\g^{3}}\right) \a^{12}  \notag\\
  &\qquad\qquad\qquad\qquad
  + \left({120}\b^{6} - {{36}\g^{2}}\b^{3} 
  + {{15}\g^{4}}\right) \a^{8} \notag\\
  &\qquad\qquad\qquad\qquad
  + \left({{-48}\g}\b^{6} + {{12}\g^{3}}\b^{3} 
  + {{6}\g^{5}}\right) \a^{4} \notag\\
  &\qquad\qquad\qquad\qquad
  + \left({64}\b^{9} + {{48}\g^{2}}\b^{6} + {{12}\g^{4}}\b^{3} 
  + {\g^{6}}\right)
   \bigr) \\
  %%%%%%%%%%%%%%%%%%%%
  x &= 3\cdot\bigl(
  {\a^{8} + {{2}\g}\a^{4} + {4}\b^{3} + {\g^{2}}}  \bigr) \\
  %%%%%%%%%%%%%%%%%%%%
  y &= {{{-108}\b^{3}}\a^{4}} \\
  %%%%%%%%%%%%%%%%%%%%
  \Disc  & = 4A^3+27B^2 \notag\\
  &= 2^8 \cdot 3^{12}\cdot \b^9 \a^{8} \bigl(
    {{\g}}\a^{12} + \left(-\b^{3} + {{3}\g^{2}}\right) \a^{8} 
   + \left({{-20}\g}\b^{3} + {{3}\g^{3}}\right) \a^{4} \notag\\
  &\qquad\qquad\qquad\qquad
   + \left({16}\b^{6} + {{8}\g^{2}}\b^{3} + {\g^{4}}\right)  \bigr)
  \label{equation:formulaforDisc}
\end{align}
Then the elliptic divisibility sequence~$[1,\a,\b,\a\g,\cdots]$ is associated
to the elliptic curve~$E$ and rational point~$P\in E$ given by the
equations
\[
  E:Y^2=X^3+AX+B\qquad\text{and}\qquad P=(x,y)\in E.
\]
Up to a substitutation and a linear shift of coordinates,
the formulas \eqref{equation:formulaforA}--\eqref{equation:formulaforDisc}
for~$A$,~$B$,~$x$,~$y$, and~$\Disc$
may be found in Ward's paper~\cite{Ward1} as, respectively,
equations~(13.6),~(13.7),~(13.5),~(13.1), and~(19.3).


\begin{thebibliography}{99}

\itemsep=\smallskipamount

\bibitem{BPvdP}
J.P. B\'ezivin, A. Peth\"o, A.J. van der Poorten, A full
characterization of divisibility sequences, \emph{Amer. J. of
Math.} \textbf{112} (1990), 985--1001.

\bibitem{Durst}
L.K. Durst, The apparition problem for equianharmonic divisibility,
Proc. Nat. Acad. Sci. U. S. A. \textbf{38} (1952), 330--333.

\bibitem{EGW}
M. Einsiedler, G. Everest, T. Ward, Primes in elliptic divisibility
sequences, \emph{LMS J. Comput. Math.} \textbf{4} (2001), 1--13,
{electronic}.

\bibitem{EPSW}
G. Everest, A. van der Poorten, I. Shparlinski, T. Ward,
\emph{Recurrence sequences}, Mathematical Surveys and Monographs
104, AMS, Providence, RI, 2003.

\bibitem{EPW}
G. Everest, Y. Puri, T. Ward,
Integer sequences counting periodic points,
\emph{Journal of Integer Sequences} (electronic) 
\textbf{5} (2002), 1-10. (ArXiv math.NT/0204173)

\bibitem{EW1}
G. Everest, T. Ward, Primes in divisibility sequences,
\emph{Cubo Mat. Educ.} \textbf{3} (2001), 245--259.

\bibitem{EW2}
\bysame
The canonical height of an algebraic point on an elliptic curve,
\emph{New York J. Math.} 6 (2000), 331-342, (electronic).

%% \bibitem{Hartshorne}
%% R. Hartshorne, \emph{Algebraic Geometry}, GTM 52, Springer-Verlag, New
%% York, 1977.

%% \bibitem{LangDG} 
%% S. Lang, \emph{Fundamentals of Diophantine Geometry},
%% Springer-Verlag, New York, 1983.

\bibitem{MT}
B. Mazur, J. Tate, The $p$-adic sigma function, \emph{Duke Math. J.}
\textbf{62} (1991), 663--688.

\bibitem{SHIP}
R. Shipsey, 
%% \emph{Elliptic divisibility sequences and elliptic curves}, 
%% preprint 2001 (unpublished).
Elliptic divisibility sequences, Ph.D. thesis, Goldsmith's College
(University of London), 2000.

\bibitem{AEC}
J.H. Silverman, \emph{The arithmetic of elliptic curves}, GTM~106, 
Springer-Ver\-lag, New York, 1986.

\bibitem{ATAEC}
\bysame
\emph{Advanced topics in the arithmetic of
elliptic curves}, GTM~151, Springer-Ver\-lag, New York, 1994.

\bibitem{SilvermanGCDinFGgps}
\bysame
Algebraic divisibility sequences and Vojta's conjecture
for exceptional divisors, in preparation.

\bibitem{SW}
C.S. Swart, 
Elliptic divisibility sequences, Ph.D. thesis, Royal Holloway
(University of London), 2003.

\bibitem{Vojta}
P. Vojta, \emph{Diophantine approximations and value distribution theory},
Lecture Notes in Mathematics~1239, Springer-Verlag, Berlin, 1987.

\bibitem{Ward1}
M. Ward, Memoir on elliptic divisibility sequences, \emph{Amer. J. Math.}
\textbf{70} (1948), 31--74.

\bibitem{Ward2}
M. Ward, The law of repetition of primes in an elliptic divisibility
sequence, \emph{Duke Math. J.} \textbf{15} (1948), 941--946.

\end{thebibliography}
\end{document}